\documentclass[preprint,12pt]{elsarticle}
\usepackage{amsmath}
\usepackage{amssymb}

\usepackage{amssymb}

\newcommand{\THM}[1]{%
  \par\vspace{5mm}\par\noindent {\scshape #1.~}}
\newcommand{\COR}[1]{%
    \par\vspace{5mm}\par\noindent {\scshape #1.}}

\usepackage[symbol]{footmisc}

\begin{document}

\begin{frontmatter}

\title{Identities of symmetry for Bernoulli polynomials
    arising from quotients of Volkenborn integrals invariant under $S_{3}$}

 \author[dsk]{Dae San Kim\corref{cor1}}
 \ead{dskim@sogang.ac.kr}

\author[khp]{ Kyoung Ho Park\corref{cor2}}
 \ead{sagamath@yahoo.co.kr}

\cortext[cor1]{Principal Corresponding author}
\cortext[cor2]{Corresponding author}

\address{Department of Mathematics, Sogang University,
 Seoul 121-742, South Korea}

\tnotetext[t1] {This work was supported by National Foundation of Korea Grant funded
by the Korean Government(2009-0072514).}

\begin{abstract}
In this paper, we derive eight basic identities of symmetry in three variables related to Bernoulli polynomials and power sums. These and most of their corollaries are new, since there have been results only about identities of symmetry in two variables. These abundance of symmetries shed new light even on the existing identities so as to yield some further interesting ones. The derivations of identities are based on the $p$-adic integral expression of the generating function for the Bernoulli polynomials and the quotient of integrals that can be expressed as the exponential generating function for the power sums.
\end{abstract}


\end{frontmatter}
\newcommand{\SECTION}[2]{%
  \vspace{5mm}\par\noindent $\S$#1. {\bf #2}~}

\SECTION{I}{Introduction and preliminaries }
\par
\vspace{3mm}
Let $p$  be a fixed prime. Throughout this paper, $\mathbb{Z}_{p},\mathbb{Q}_{p},\mathbb{C}_{p}$ will respectively denote the ring of $p$-adic integers, the field of $p$-adic rational numbers and the completion of the algebraic closure of $\mathbb{Q}_{p}$. For a uniformly differentiable (also called continuously differentiable) function  $f: \mathbb{Z}_{p} \rightarrow \mathbb{C}_{p}$ (cf. [4]), the Volkenborn integral of $f$ is defined by
 $$
 \int_{\mathbb{Z}_{p}}f(z)d\mu(z)=\lim_{N \rightarrow \infty } \frac{1}{p^{N}}\sum_{j=0}^{p^{N}-1}
f(j).
$$

Then it is easy to see that
\begin{equation}
 \int_{\mathbb{Z}_{p}}f(z+1)d\mu(z)=\int_{\mathbb{Z}_{p}}f(z)d\mu(z)+f^{\prime}(0).
\end{equation}
Let $|\cdot |_{p}$ be the normalized absolute value of $\mathbb{C}_{p}$, such that $|p|_{p}=\frac{1}{p}$, and let
\begin{equation}
E=\{ t \in \mathbb{C}_{p}| \ |t|_{p} < p^{-\frac{1}{p-1}}\}.
\end{equation}

Then, for each fixed $t \in E$, the function $f(z)=e^{zt}$ is analytic on $\mathbb{Z}_{p}$ and by applying (1) to this  $f$, we get the $p$-adic integral expression of the generating function for Bernoulli numbers $B_{n}:$

\begin{equation}
 \int_{\mathbb{Z}_{p}}e^{zt}d\mu(z)=\frac{t}{e^{t}-1}=\sum_{n=0}^{\infty}B_{n}\frac{t^{n}}{n!} \quad (t \in E).
\end{equation}

So we have the following $p$-adic integral expression of the generating function for the Bernoulli polynomials $B_{n}(x):$
\begin{equation}
 \int_{\mathbb{Z}_{p}}e^{(x+z)t}d\mu(z)=\frac{t}{e^{t}-1}e^{xt}=\sum_{n=0}^{\infty}B_{n}(x)\frac{t^{n}}{n!} \quad (t \in E, x\in \mathbb{Z}_{p}).
\end{equation}

Here and throughout this paper, we will have many instances to be able to interchange integral and infinite sum. That is justified by Proposition 55.4 in [4].
Let $S_{k}(n)$ denote the $k$-th power sum of the first $n+1$ nonnegative integers, namely
\begin{equation}
S_{k}(n)=\sum_{i=0}^{n}i^{k}=0^{k}+1^{k}+\cdots+n^{k}.
\end{equation}

In particular,
\begin{equation*}
\tag{6}
S_{0}(n)=n+1,
\quad
S_{k}(0)=
\begin{cases}
1, & \rm {for}\ \it{k}=\rm{0},\\
0, & \rm {for}\ \it{k}>\rm{0}.
\end{cases}
\end{equation*}

From (3) and (5), one easily derives the following identities: for $w \in \mathbb{Z}_{>0}$,
\begin{equation}
\tag{7}
\frac{w \int_{\mathbb{Z}_{p}}e^{xt}d\mu(x)}{\int_{\mathbb{Z}_{p}}e^{ w yt}d\mu(y)}=\sum_{i=0}^{w -1}e^{it}=\sum_{k=0}^{\infty}S_{k}(w -1)\frac{t^{k}}{k!}\quad ( t \in E ).
\end{equation}
In what follows, we will always assume that the Volkenborn integrals of the various exponential functions on $\mathbb{Z}_{p}$ are defined for $ t\in E$ (cf.(2)), and therefore it will not be mentioned.
\par
[1]-[3], [5] and [6] are some of the previous works on identities of symmetry in two variables involving Bernoulli polynomials and power sums. For the brief history, one is referred to those papers.
\par
In this paper, we will produce  8 basic identities of symmetry in three variables $w_{1}, w_{2}, w_{3}$ related to Bernoulli polynomials and power sums(cf. (44), (45), (48), (51), (55), (57), (59), (60)). These and most of their corollaries seem to be new,  since there have been results only about identities of symmetry in two variables in the literature. These abundance of symmetries shed new light even on the existing identities. For instance, it has been known that (8) and (9) are equal(cf. [6, Cor.1], [3, Cor.2]) and (10) and (11) are so (cf. [3, (13)], [6, Cor.4]).
In fact, (8)-(11) are all equal, as they can be derived from one and the same $p$-adic integral. Perhaps, this was
neglected to mention in [3]. Also, we have a bunch of new identities in (12)-(15). All of these were obtained
as corollaries (cf. Cor.9, 12, 15) to some of the basic identities by specializing the variable $w_{3}$ as 1.
Those would not be unearthed if more symmetries had not been available.

\begin{align}
\tag{8}
&\sum_{k=0}^{n}\binom{n}{k}B_{k}(w_{1}y_{1})S_{n-k}(w_{2}-1)w_{1}^{n-k}w_{2}^{k-1}\\
\tag{9}
&= \sum_{k=0}^{n}\binom{n}{k}B_{k}(w_{2}y_{1})S_{n-k}(w_{1}-1)w_{2}^{n-k}w_{1}^{k-1}\\
\tag{10}
&= w_{1}^{n-1}\sum_{i=0}^{w_{1}-1}B_{n}(w_{2}y_{1}+\frac{w_{2}}{w_{1}}i)\\
\tag{11}
&= w_{2}^{n-1}\sum_{i=0}^{w_{2}-1}B_{n}(w_{1}y_{1}+\frac{w_{1}}{w_{2}}i)\\
\tag{12}
&=\sum_{k+\ell+m=n}\binom{n}{k,\ell,m}B_{k}(y_{1})S_{\ell}(w_{1}-1)S_{m}(w_{2}-1)w_{1}^{k+m-1}w_{2}^{k+\ell-1}
\end{align}
\begin{align}
\tag{13}
&= w_{1}^{n-1}\sum_{k=0}^{n}\binom{n}{k}\sum_{i=0}^{w_{1}-1}B_{k}(y_{1}+\frac{i}{w_{1}})S_{n-k}(w_{2}-1)w_{2}^{k-1}\\
\tag{14}
&= w_{2}^{n-1}\sum_{k=0}^{n}\binom{n}{k}\sum_{i=0}^{w_{2}-1}B_{k}(y_{1}+\frac{i}{w_{2}})S_{n-k}(w_{1}-1)w_{1}^{k-1}\\
\tag{15}
&= (w_{1}w_{2})^{n-1}\sum_{i=0}^{w_{1}-1}\sum_{j=0}^{w_{2}-1}B_{n}(y_{1}+\frac{i}{w_{1}}+\frac{j}{w_{2}}).
\end{align}

The derivations of identities is based on the $p$-adic integral expression of the generating function for the Bernoulli polynomials in (4) and the quotient of integrals in (7) that can be expressed as the exponential generating function for the power sums. We indebted this idea to the paper [3].

\SECTION{II}{Several types of quotients of Volkenborn integrals}
\par
\vspace{3mm}
Here we will introduce several types of quotients of Volkenborn integrals on $\mathbb{Z}_{p}$ or $\mathbb{Z}_{p}^{3}$ from which some interesting identities follow owing to the built-in symmetries in $w_{1},w_{2},w_{3}$. In the following, $w_{1},w_{2},w_{3}$ are positive integers and all of the explicit expressions of integrals in (17), (19), (21), and (23) are obtained from the identity in (3).
\vspace{5mm}
\par
$(a)$ Type $\Lambda_{23}^{i}$ (for  $i=0,1,2,3$ )
\vspace{3mm}
\begin{align}
\tag{16}
I(\Lambda_{23}^{i})
&=\frac{\int_{\mathbb{Z}_{p}^{3}}
e^{(w_{2}w_{3}x_{1}+w_{1}w_{3}x_{2}+w_{1}w_{2}x_{3}+w_{1}w_{2}w_{3}(\sum_{j=1}^{3-i}y_{j}))t}
d\mu(x_{1})d\mu(x_{2})d\mu(x_{3}) }{(\int_{\mathbb{Z}_{p}}e^{w_{1}w_{2}w_{3}x_{4}t}d\mu(x_{4}))^{i}} \\
\tag{17}
&= \frac{(w_{1}w_{2}w_{3})^{2-i}t^{3-i}e^{w_{1}w_{2}w_{3}(\sum_{j=1}^{3-i}y_{j})t}(e^{w_{1}w_{2}w_{3}t}-1)^{i}}
{(e^{w_{2}w_{3}t}-1)(e^{w_{1}w_{3}t}-1)(e^{w_{1}w_{2}t}-1)};
\end{align}
\vspace{3mm}

\par
$(b)$ Type $\Lambda_{13}^{i}$ (for  $i=0,1,2,3$ )

\begin{align}
\tag{18}
I(\Lambda_{13}^{i})
&=\frac{\int_{\mathbb{Z}_{p}^{3}}
e^{(w_{1}x_{1}+w_{2}x_{2}+w_{3}x_{3}+w_{1}w_{2}w_{3}(\sum_{j=1}^{3-i}y_{j}))t}
d\mu(x_{1})d\mu(x_{2})d\mu(x_{3}) }{(\int_{\mathbb{Z}_{p}}e^{w_{1}w_{2}w_{3}x_{4}t}d\mu(x_{4}))^{i}}
\end{align}
\begin{align}
\tag{19}
&= \frac{(w_{1}w_{2}w_{3})^{1-i}t^{3-i}e^{w_{1}w_{2}w_{3}(\sum_{j=1}^{3-i}y_{j})t}(e^{w_{1}w_{2}w_{3}t}-1)^{i}}
{(e^{w_{1}t}-1)(e^{w_{2}t}-1)(e^{w_{3}t}-1)};
\end{align}

\vspace{3mm}

\par

$(c-0)$ Type $\Lambda_{12}^{0}$

\begin{align}
\tag{20}
I(\Lambda_{12}^{0})
&=\int_{\mathbb{Z}_{p}^{3}}
e^{(w_{1}x_{1}+w_{2}x_{2}+w_{3}x_{3}+w_{2}w_{3}y+w_{1}w_{3}y+w_{1}w_{2}y)t} d\mu(x_{1})d\mu(x_{2})d\mu(x_{3}) \\
\tag{21}
&= \frac{w_{1}w_{2}w_{3}t^{3}e^{(w_{2}w_{3}+w_{1}w_{3}+w_{1}w_{2})yt}}
{(e^{w_{1}t}-1)(e^{w_{2}t}-1)(e^{w_{3}t}-1)};
\end{align}
\vspace{3mm}

\par

$(c-1)$ Type $\Lambda_{12}^{1}$

\begin{align}
\tag{22}
I(\Lambda_{12}^{1})
&=\frac{\int_{\mathbb{Z}_{p}^{3}}
e^{(w_{1}x_{1}+w_{2}x_{2}+w_{3}x_{3})t}
d\mu(x_{1})d\mu(x_{2})d\mu(x_{3}) }{\int_{\mathbb{Z}_{p}^{3}}
e^{(w_{2}w_{3}z_{1}+w_{1}w_{3}z_{2}+w_{1}w_{2}z_{3})t}
d\mu(z_{1})d\mu(z_{2})d\mu(z_{3})} \\
\tag{23}
&= \frac{(w_{1}w_{2}w_{3})^{-1}(e^{w_{2}w_{3}t}-1)(e^{w_{1}w_{3}t}-1)(e^{w_{1}w_{2}t}-1)}
{(e^{w_{1}t}-1)(e^{w_{2}t}-1)(e^{w_{3}t}-1)}.
\end{align}
\par
All of the above $p$-adic integrals of various types are invariant under all permutations of  $w_{1},w_{2},w_{3}$, as one can see either from $p$-adic integral representations in (16), (18), (20), and (22) or from their explicit evaluations in (17), (19), (21), and (23).

\SECTION{III}{Identities for Bernoulli polynomials}
\par
\vspace{3mm}

$(a-0)$ First, let's consider Type $\Lambda_{23}^{i}$, for each $i=0,1,2,3$. The following results can be easily obtained from (4) and (7).
\begin{align}
\notag
&I(\Lambda_{23}^{0})\\
\notag
&=\int_{\mathbb{Z}_{p}}e^{w_{2}w_{3}(x_{1}+w_{1}y_{1})t}d\mu(x_{1})
\int_{\mathbb{Z}_{p}}e^{w_{1}w_{3}(x_{2}+w_{2}y_{2})t}d\mu(x_{2})
\int_{\mathbb{Z}_{p}}e^{w_{1}w_{2}(x_{3}+w_{3}y_{3})t}d\mu(x_{3})\\
\tag{24}
&=(\sum_{k=0}^{\infty}\frac{B_{k}(w_{1}y_{1})}{k!}(w_{2}w_{3}t)^{k})
(\sum_{\ell=0}^{\infty}\frac{B_{\ell}(w_{2}y_{2})}{\ell!}(w_{1}w_{3}t)^{\ell})
(\sum_{m=0}^{\infty}\frac{B_{m}(w_{3}y_{3})}{m!}(w_{1}w_{2}t)^{m})
\end{align}

\begin{equation*}
=\sum_{n=0}^{\infty}\Big( \sum_{k+\ell+m=n}\binom{n}{k,\ell,m}B_{k}(w_{1}y_{1})B_{\ell}(w_{2}y_{2})B_{m}(w_{3}y_{3})
w_{1}^{\ell+m}w_{2}^{k+m}w_{3}^{k+\ell}\Big )\frac{t^{n}} {n!},
\end{equation*}
where the inner sum is over all nonnegative integers $k, \ell, m$,with $k+\ell+m=n$, and
\begin{equation*}
\tag{25}
\binom{n}{k,\ell,m}=\frac{n!}{k!\ell!m!}.
\end{equation*}
\vspace{5mm}
\par
$(a-1)$  Here we write $I(\Lambda_{23}^{1})$ in two different ways:

\begin{align}
\notag
(1) \ I(\Lambda_{23}^{1})
&=\frac{1}{w_{3}}\int_{\mathbb{Z}_{p}}e^{w_{2}w_{3}(x_{1}+w_{1}y_{1})t}d\mu(x_{1})
\int_{\mathbb{Z}_{p}}e^{w_{1}w_{3}(x_{2}+w_{2}y_{2})t}d\mu(x_{2})\\
\tag{26}
&\times  \frac{w_{3}\int_{\mathbb{Z}_{p}}e^{w_{1}w_{2}x_{3}t}d\mu(x_{3})}
 {\int_{\mathbb{Z}_{p}}e^{w_{1}w_{2}w_{3}x_{4}t}d\mu(x_{4})}\\
\notag
&=\frac{1}{w_{3}}\Big(\sum_{k=0}^{\infty}B_{k}(w_{1}y_{1})\frac{(w_{2}w_{3}t)^{k}}{k!}\Big)
\Big(\sum_{\ell=0}^{\infty}B_{\ell}(w_{2}y_{2})\frac{(w_{1}w_{3}t)^{\ell}}{\ell!}\Big)\\
\notag
&\times \Big(\sum_{m=0}^{\infty}S_{m}(w_{3}-1)\frac{(w_{1}w_{2}t)^{m}}{m!}\Big)\\
\notag
&=\sum_{n=0}^{\infty}\Big( \sum_{k+\ell+m=n}\binom{n}{k,\ell,m}B_{k}(w_{1}y_{1})B_{\ell}(w_{2}y_{2})S_{m}(w_{3}-1)\\
\tag{27}
&\times w_{1}^{\ell+m}w_{2}^{k+m}w_{3}^{k+\ell -1}\Big )\frac{t^{n}} {n!}.\\
\notag
\end{align}

(2) Invoking (7), (26) can also be written as

\begin{align}
\notag
&I(\Lambda_{23}^{1})\\
\notag
&=\frac{1}{w_{3}}\sum_{i=0}^{w_{3}-1}\int_{\mathbb{Z}_{p}}e^{w_{2}w_{3}(x_{1}+w_{1}y_{1})t}d\mu(x_{1})
\int_{\mathbb{Z}_{p}}e^{w_{1}w_{3}(x_{2}+w_{2}y_{2}+\frac{w_{2}}{w_{3}}i)t}d\mu(x_{2})\\
\notag
\end{align}
\begin{align}
\notag
&=\frac{1}{w_{3}}\sum_{i=0}^{w_{3}-1}\Big(\sum_{k=0}^{\infty}B_{k}(w_{1}y_{1})\frac{(w_{2}w_{3}t)^{k}}{k!}\Big)
\Big(\sum_{\ell=0}^{\infty}B_{\ell}(w_{2}y_{2}+\frac{w_{2}}{w_{3}}i)\frac{(w_{1}w_{3}t)^{\ell}}{\ell!}\Big)\\
\tag{28}
&=\sum_{n=0}^{\infty}\Big( w_{3}^{n-1}\sum_{k=0}^{n}\binom{n}{k}B_{k}(w_{1}y_{1})
\sum_{i=0}^{w_{3}-1}B_{n-k}(w_{2}y_{2}+\frac{w_{2}}{w_{3}}i)w_{1}^{n-k}w_{2}^{k}\Big)\frac{t^{n}} {n!}.\\
\notag
\end{align}


\par

$(a-2)$  Here we write $I(\Lambda_{23}^{2})$ in three different ways:

\begin{align}
\notag
(1)  \ I(\Lambda_{23}^{2})
&=\frac{1}{w_{2}w_{3}}\int_{\mathbb{Z}_{p}}e^{w_{2}w_{3}(x_{1}+w_{1}y_{1})t}d\mu(x_{1})
\times
\frac{w_{2}\int_{\mathbb{Z}_{p}}e^{w_{1}w_{3}x_{2}t}d\mu(x_{2}) }
{\int_{\mathbb{Z}_{p}}e^{w_{1}w_{2}w_{3}x_{4}t}d\mu(x_{4}) }\\
\tag{29}
& \times \frac{w_{3}\int_{\mathbb{Z}_{p}}e^{w_{1}w_{2}x_{3}t}d\mu(x_{3}) }
{\int_{\mathbb{Z}_{p}}e^{w_{1}w_{2}w_{3}x_{4}t}d\mu(x_{4}) }\\
\notag
&=\frac{1}{w_{2}w_{3}}\Big(\sum_{k=0}^{\infty}B_{k}(w_{1}y_{1})\frac{(w_{2}w_{3}t)^{k}}{k!}\Big)
\Big(\sum_{\ell=0}^{\infty}S_{\ell}(w_{2}-1)\frac{(w_{1}w_{3}t)^{\ell}}{\ell!}\Big)\\
\notag
& \times \Big(\sum_{m=0}^{\infty}S_{m}(w_{3}-1)\frac{(w_{1}w_{2}t)^{m}}{m!}\Big)\\
\notag
&=\sum_{n=0}^{\infty}\Big( \sum_{k+\ell+m=n}\binom{n}{k,\ell,m}B_{k}(w_{1}y_{1})S_{\ell}(w_{2}-1)S_{m}(w_{3}-1)\\
\tag{30}
& \times w_{1}^{\ell+m}w_{2}^{k+m-1}w_{3}^{k+\ell -1} \Big )\frac{t^{n}} {n!}.\\
\notag
\end{align}

(2) Invoking (7), (29) can also be written as
\begin{align}
\notag
&I(\Lambda_{23}^{2})\\
\tag{31}
&=\frac{1}{w_{2}w_{3}}\sum_{i=0}^{w_{2}-1}\int_{\mathbb{Z}_{p}}e^{w_{2}w_{3}(x_{1}+w_{1}y_{1}+\frac{w_{1}}{w_{2}}i)t}d\mu(x_{1})
\times
\frac{w_{3}\int_{\mathbb{Z}_{p}}e^{w_{1}w_{2}x_{3}t}d\mu(x_{3}) }
{\int_{\mathbb{Z}_{p}}e^{w_{1}w_{2}w_{3}x_{4}t}d\mu(x_{4}) }\\
\notag
&=\frac{1}{w_{2}w_{3}}\sum_{i=0}^{w_{2}-1}\Big(\sum_{k=0}^{\infty}B_{k}(w_{1}y_{1}+\frac{w_{1}}{w_{2}}i)\frac{(w_{2}w_{3}t)^{k}}{k!}\Big)
\Big(\sum_{\ell=0}^{\infty}S_{\ell}(w_{3}-1)\frac{(w_{1}w_{2}t)^{\ell}}{\ell!}\Big)\\
\tag{32}
&=\sum_{n=0}^{\infty}\Big(w_{2}^{n-1} \sum_{k=0}^{n}\binom{n}{k}
\sum_{i=0}^{w_{2}-1}B_{k}(w_{1}y_{1}+\frac{w_{1}}{w_{2}}i)S_{n-k}(w_{3}-1)
w_{1}^{n-k}w_{3}^{k-1}\Big )\frac{t^{n}} {n!}.
\end{align}

(3) Invoking (7) once again, (31) can be written as
\begin{align}
\notag
&I(\Lambda_{23}^{2})\\
\notag
&=\frac{1}{w_{2}w_{3}}\sum_{i=0}^{w_{2}-1}\sum_{j=0}^{w_{3}-1}
\int_{\mathbb{Z}_{p}}e^{w_{2}w_{3}(x_{1}+w_{1}y_{1}+\frac{w_{1}}{w_{2}}i+\frac{w_{1}}{w_{3}}j)t}d\mu(x_{1})\\
\notag
&=\frac{1}{w_{2}w_{3}}\sum_{i=0}^{w_{2}-1} \sum_{j=0}^{w_{3}-1}  \Big(\sum_{n=0}^{\infty} B_{n}(w_{1}y_{1}+\frac{w_{1}}{w_{2}}i+\frac{w_{1}}{w_{3}}j)\frac{(w_{2}w_{3}t)^{n}}{n!}\Big)\\
\tag{33}
&=\sum_{n=0}^{\infty}\Big((w_{2}w_{3})^{n-1} \sum_{i=0}^{w_{2}-1}
\sum_{j=0}^{w_{3}-1}B_{n}(w_{1}y_{1}+\frac{w_{1}}{w_{2}}i+\frac{w_{1}}{w_{3}}j )\Big )\frac{t^{n}} {n!}.
\end{align}

\vspace{3mm}
$(a-3)$
\begin{align}
\notag
I(\Lambda_{23}^{3})
&=\frac{1}{w_{1}w_{2}w_{3}}
\times
\frac{w_{1}\int_{\mathbb{Z}_{p}}e^{w_{2}w_{3}x_{1}t}d\mu(x_{1}) }
{\int_{\mathbb{Z}_{p}}e^{w_{1}w_{2}w_{3}x_{4}t}d\mu(x_{4}) }
\times
\frac{w_{2}\int_{\mathbb{Z}_{p}}e^{w_{1}w_{3}x_{2}t}d\mu(x_{2}) }
{\int_{\mathbb{Z}_{p}}e^{w_{1}w_{2}w_{3}x_{4}t}d\mu(x_{4}) }\\
\notag
&\times
\frac{w_{3}\int_{\mathbb{Z}_{p}}e^{w_{1}w_{2}x_{3}t}d\mu(x_{3}) }
{\int_{\mathbb{Z}_{p}}e^{w_{1}w_{2}w_{3}x_{4}t}d\mu(x_{4}) } \\
\notag
&=\frac{1}{w_{1}w_{2}w_{3}}\Big(\sum_{k=0}^{\infty}S_{k}(w_{1}-1)\frac{(w_{2}w_{3}t)^{k}}{k!}\Big)
\Big(\sum_{\ell=0}^{\infty}S_{\ell}(w_{2}-1)\frac{(w_{1}w_{3}t)^{\ell}}{\ell!}\Big)\\
\notag
&\times \Big(\sum_{m=0}^{\infty}S_{m}(w_{3}-1)\frac{(w_{1}w_{2}t)^{m}}{m!} \Big)\\
\notag
&=\sum_{n=0}^{\infty}\Big( \sum_{k+\ell+m=n}\binom{n}{k,\ell,m}S_{k}(w_{1}-1)S_{\ell}(w_{2}-1)S_{m}(w_{3}-1)\\
\tag{34}
& \times w_{1}^{\ell+m-1}w_{2}^{k+m-1}w_{3}^{k+\ell -1}\Big )\frac{t^{n}} {n!}.\\
\notag
\end{align}

\vspace{3mm}
\par
$(b)$ For Type $\Lambda_{13}^{i} (i= 0, 1,2,3)$, we may consider the analogous things to the ones in
$(a-0),(a-1),(a-2)$, and $(a-3)$. However, these do not lead us to new identities. Indeed, if we substitute  $w_{2}w_{3},w_{1}w_{3},w_{1}w_{2}$ respectively for $w_{1},w_{2},w_{3}$ in (16), this amounts to replacing $t$ by
$w_{1}w_{2}w_{3}t$ in (18). So, upon replacing $w_{1},w_{2},w_{3}$  respectively by  $w_{2}w_{3},w_{1}w_{3},w_{1}w_{2}$ and dividing by $(w_{1}w_{2}w_{3})^{n}$, in each of the expressions of (24), (27), (28), (30), (32)-(34), we will get the corresponding symmetric identities for  Type $\Lambda_{13}^{i}(i=0,1,2,3)$.
\vspace{5mm}
\par
$(c-0)$
\begin{equation*}
\tag{35}
\begin{split}
I(\Lambda_{12}^{0})
&=\int_{\mathbb{Z}_{p}}e^{w_{1}(x_{1}+w_{2}y)t}d\mu(x_{1})
\int_{\mathbb{Z}_{p}}e^{w_{2}(x_{2}+w_{3}y)t}d\mu(x_{2})
\int_{\mathbb{Z}_{p}}e^{w_{3}(x_{3}+w_{1}y)t}d\mu(x_{3})\\
&=\Big(\sum_{k=0}^{\infty}\frac{B_{k}(w_{2}y)}{k!}(w_{1}t)^{k}\Big)
\Big(\sum_{\ell=0}^{\infty}\frac{B_{\ell}(w_{3}y)}{\ell!}(w_{2}t)^{\ell}\Big)
\Big(\sum_{m=0}^{\infty}\frac{B_{m}(w_{1}y)}{m!}(w_{3}t)^{m}\Big)\\
&=\sum_{n=0}^{\infty}\Big( \sum_{k+\ell+m=n}\binom{n}{k,\ell,m}B_{k}(w_{2}y)B_{\ell}(w_{3}y)B_{m}(w_{1}y)
w_{1}^{k}w_{2}^{\ell}w_{3}^{m}\Big )\frac{t^{n}} {n!},\\
\end{split}
\end{equation*}

\vspace{5mm}

$(c-1)$

\begin{align}
\notag
& I(\Lambda_{12}^{1}) \\
\notag
&=\frac{1}{w_{1}w_{2}w_{3}}
\frac{w_{2}\int_{\mathbb{Z}_{p}}e^{w_{1}x_{1}t}d\mu(x_{1}) }
{\int_{\mathbb{Z}_{p}}e^{w_{1}w_{2}z_{3}t}d\mu(z_{3}) }
\times
\frac{w_{3}\int_{\mathbb{Z}_{p}}e^{w_{2}x_{2}t}d\mu(x_{2}) }
{\int_{\mathbb{Z}_{p}}e^{w_{2}w_{3}z_{1}t}d\mu(z_{1}) }
\times
\frac{w_{1}\int_{\mathbb{Z}_{p}}e^{w_{3}x_{3}t}d\mu(x_{3}) }
{\int_{\mathbb{Z}_{p}}e^{w_{3}w_{1}z_{2}t}d\mu(z_{2}) } \\
\notag
&=\frac{1}{w_{1}w_{2}w_{3}}\Big(\sum_{k=0}^{\infty}S_{k}(w_{2}-1)\frac{(w_{1}t)^{k}}{k!}\Big)
\Big(\sum_{\ell=0}^{\infty}S_{\ell}(w_{3}-1)\frac{(w_{2}t)^{\ell}}{\ell!}\Big) \\
\notag
& \times \Big(\sum_{m=0}^{\infty}S_{m}(w_{1}-1)\frac{(w_{3}t)^{m}}{m!}\Big) \\
\tag{36}
& =\sum_{n=0}^{\infty}\Big( \sum_{k+\ell+m=n}\binom{n}{k,\ell,m}S_{k}(w_{2}-1)S_{\ell}(w_{3}-1)S_{m}(w_{1}-1)\\
\notag
& \times w_{1}^{k-1}w_{2}^{\ell-1}w_{3}^{m -1}\Big )\frac{t^{n}} {n!}.
\notag
\end{align}

\SECTION{IV}{Main theorems}
\par
\vspace{3mm}
As we noted earlier in the last paragraph of Section II, the various types of quotients of Volkenborn integrals are invariant under any permutation of $w_{1}, w_{2},w_{3}$. So the corresponding expressions in Section III are also invariant under any permutation of $w_{1}, w_{2},w_{3}$. Thus our results about identities of symmetry will be immediate consequences of this observation.
\par
However, not all permutations of an expression in Section III yield distinct ones. In fact, as these expressions are obtained by permuting $w_{1}, w_{2},w_{3}$ in a single one labeled by them, they can be viewed as a group in a natural manner and hence it is isomorphic to a quotient of $S_{3}$. In particular, the number of possible distinct expressions are 1,2,3, or 6.(a-0),(a-1(1)),(a-1(2)), and (a-2(2)) give the full six identities of symmetry, (a-2(1)) and (a-2(3)) yield three identities of symmetry, and (c-0) and (c-1) give two identities of symmetry, while the expression in (a-3) yields no identities of symmetry.
\par
Here we will just consider the cases of Theorems 8 and 17, leaving the others as easy exercises for the reader. As for the case of Theorem 8, in addition to (50)-(52), we get the following three ones:
\begin{equation*}
\tag{37}
\sum_{k+\ell+m=n}\binom{n}{k,\ell,m}B_{k}(w_{1}y_{1})S_{\ell}(w_{3}-1)S_{m}(w_{2}-1)
w_{1}^{\ell+m}w_{3}^{k+m-1}w_{2}^{k+\ell -1},
\end{equation*}

\begin{equation*}
\tag{38}
\sum_{k+\ell+m=n}\binom{n}{k,\ell,m}B_{k}(w_{2}y_{1})S_{\ell}(w_{1}-1)S_{m}(w_{3}-1)
w_{2}^{\ell+m}w_{1}^{k+m-1}w_{3}^{k+\ell -1},
\end{equation*}

\begin{equation*}
\tag{39}
\sum_{k+\ell+m=n}\binom{n}{k,\ell,m}B_{k}(w_{3}y_{1})S_{\ell}(w_{2}-1)S_{m}(w_{1}-1)
w_{3}^{\ell+m}w_{2}^{k+m-1}w_{1}^{k+\ell -1}.
\end{equation*}

But, by interchanging $\ell$ and $m$, we see that (37), (38), and (39) are respectively equal to (50), (51), and (52).\\
As to Theorem 17, in addition to (60) and (61), we have:

\begin{align}
\tag{40}
\sum_{k+\ell+m=n}\binom{n}{k,\ell,m}S_{k}(w_{2}-1)S_{\ell}(w_{3}-1)S_{m}(w_{1}-1)
w_{1}^{k-1}w_{2}^{\ell-1}w_{3}^{m-1},\\
\tag{41}
\sum_{k+\ell+m=n}\binom{n}{k,\ell,m}S_{k}(w_{3}-1)S_{\ell}(w_{1}-1)S_{m}(w_{2}-1)
w_{2}^{k-1}w_{3}^{\ell-1}w_{1}^{m-1},\\
\tag{42}
\sum_{k+\ell+m=n}\binom{n}{k,\ell,m}S_{k}(w_{3}-1)S_{\ell}(w_{2}-1)S_{m}(w_{1}-1)
w_{1}^{k-1}w_{3}^{\ell-1}w_{2}^{m -1},\\
\tag{43}
\sum_{k+\ell+m=n}\binom{n}{k,\ell,m}S_{k}(w_{2}-1)S_{\ell}(w_{1}-1)S_{m}(w_{3}-1)
w_{3}^{k-1}w_{2}^{\ell-1}w_{1}^{m -1}.
\end{align}

However, (40) and (41) are equal to (60), as we can see by applying the permutations
$k \rightarrow \ell, \ell \rightarrow m, m \rightarrow k$ for (40) and
$k \rightarrow m, \ell \rightarrow k, m \rightarrow \ell$ for (41).
Similarly, we see that (42) and (43) are equal to (61), by applying permutations
$k \rightarrow \ell, \ell \rightarrow m, m \rightarrow k$ for (42)
and $k \rightarrow m, \ell \rightarrow k, m \rightarrow \ell$ for (43).

{\THM {Theorem {1}}}{\it Let $w_{1}, w_{2},w_{3} $ be any positive integers.
Then the following expression  is invariant under any permutation of
$w_{1}, w_{2},w_{3}$, so that it gives us six symmetries.
\rm}

\begin{equation*}
\tag{44}
\begin{split}
& \sum_{k+\ell+m=n}\binom{n}{k,\ell,m}B_{k}(w_{1}y_{1})B_{\ell}(w_{2}y_{2})B_{m}(w_{3}y_{3})
w_{1}^{\ell+m}w_{2}^{k+m}w_{3}^{k+\ell}\\
&=\sum_{k+\ell+m=n}\binom{n}{k,\ell,m}B_{k}(w_{1}y_{1})B_{\ell}(w_{3}y_{2})B_{m}(w_{2}y_{3})
w_{1}^{\ell+m}w_{3}^{k+m}w_{2}^{k+\ell}\\
&=\sum_{k+\ell+m=n}\binom{n}{k,\ell,m}B_{k}(w_{2}y_{1})B_{\ell}(w_{1}y_{2})B_{m}(w_{3}y_{3})
w_{2}^{\ell+m}w_{1}^{k+m}w_{3}^{k+\ell}\\
&=\sum_{k+\ell+m=n}\binom{n}{k,\ell,m}B_{k}(w_{2}y_{1})B_{\ell}(w_{3}y_{2})B_{m}(w_{1}y_{3})
w_{2}^{\ell+m}w_{3}^{k+m}w_{1}^{k+\ell} \\
&= \sum_{k+\ell+m=n}\binom{n}{k,\ell,m}B_{k}(w_{3}y_{1})B_{\ell}(w_{1}y_{2})B_{m}(w_{2}y_{3})
w_{3}^{\ell+m}w_{1}^{k+m}w_{2}^{k+\ell}\\
&= \sum_{k+\ell+m=n}\binom{n}{k,\ell,m}B_{k}(w_{3}y_{1})B_{\ell}(w_{2}y_{2})B_{m}(w_{1}y_{3})
w_{3}^{\ell+m}w_{2}^{k+m}w_{1}^{k+\ell}.\\
\end{split}
\end{equation*}

{\THM {Theorem {2}}}{\it Let $w_{1}, w_{2},w_{3} $ be any positive integers.
Then the following expression  is invariant under any permutation of
$w_{1}, w_{2},w_{3}$, so that it gives us six symmetries.
\rm}
\begin{equation*}
\tag{45}
\begin{split}
& \sum_{k+\ell+m=n}\binom{n}{k,\ell,m}B_{k}(w_{1}y_{1})B_{\ell}(w_{2}y_{2})S_{m}(w_{3}-1)
w_{1}^{\ell+m}w_{2}^{k+m}w_{3}^{k+\ell-1}\\
&=\sum_{k+\ell+m=n}\binom{n}{k,\ell,m}B_{k}(w_{1}y_{1})B_{\ell}(w_{3}y_{2})S_{m}(w_{2}-1)
w_{1}^{\ell+m}w_{3}^{k+m}w_{2}^{k+\ell-1}\\
&=\sum_{k+\ell+m=n}\binom{n}{k,\ell,m}B_{k}(w_{2}y_{1})B_{\ell}(w_{1}y_{2})S_{m}(w_{3}-1)
w_{2}^{\ell+m}w_{1}^{k+m}w_{3}^{k+\ell-1}\\
\end{split}
\end{equation*}

\begin{equation*}
\begin{split}
&=\sum_{k+\ell+m=n}\binom{n}{k,\ell,m}B_{k}(w_{2}y_{1})B_{\ell}(w_{3}y_{2})S_{m}(w_{1}-1)
w_{2}^{\ell+m}w_{3}^{k+m}w_{1}^{k+\ell-1} \\
&= \sum_{k+\ell+m=n}\binom{n}{k,\ell,m}B_{k}(w_{3}y_{1})B_{\ell}(w_{2}y_{2})S_{m}(w_{1}-1)
w_{3}^{\ell+m}w_{2}^{k+m}w_{1}^{k+\ell-1}\\
&= \sum_{k+\ell+m=n}\binom{n}{k,\ell,m}B_{k}(w_{3}y_{1})B_{\ell}(w_{1}y_{2})S_{m}(w_{2}-1)
w_{3}^{\ell+m}w_{1}^{k+m}w_{2}^{k+\ell-1}.\\
\end{split}
\end{equation*}

Putting $w_{3}=1$ in (45), we get the following corollary.

{\COR {Corollary {3}}} {\it Let $w_{1}, w_{2}$ be any positive integers. \rm}

\begin{equation*}
\tag{46}
\begin{split}
& \sum_{k=0}^{n}\binom{n}{k}B_{k}(w_{1}y_{1})B_{n-k}(w_{2}y_{2})w_{1}^{n-k}w_{2}^{k}\\
&=\sum_{k=0}^{n}\binom{n}{k}B_{k}(w_{2}y_{1})B_{n-k}(w_{1}y_{2})w_{2}^{n-k}w_{1}^{k}\\
&=\sum_{k+\ell+m=n}\binom{n}{k,\ell,m}B_{k}(y_{1})B_{\ell}(w_{2}y_{2})S_{m}(w_{1}-1)
w_{2}^{k+m}w_{1}^{k+\ell -1}\\
&=\sum_{k+\ell+m=n}\binom{n}{k,\ell,m}B_{k}(w_{2}y_{1})B_{\ell}(y_{2})S_{m}(w_{1}-1)
w_{2}^{\ell+m}w_{1}^{k+\ell-1} \\
&= \sum_{k+\ell+m=n}\binom{n}{k,\ell,m}B_{k}(y_{1})B_{\ell}(w_{1}y_{2})S_{m}(w_{2}-1)
w_{1}^{k+m}w_{2}^{k+\ell-1}\\
\end{split}
\end{equation*}
\begin{equation*}
= \sum_{k+\ell+m=n}\binom{n}{k,\ell,m}B_{k}(w_{1}y_{1})B_{\ell}(y_{2})S_{m}(w_{2}-1)
w_{1}^{\ell+m}w_{2}^{k+\ell-1}.
\end{equation*}

Letting further  $w_{2}=1$ in (46), we have the following corollary.

{\COR {Corollary {4}}} {\it Let $w_{1}$ be any positive integer. \rm}

\begin{equation*}
\tag{47}
\begin{split}
& \sum_{k=0}^{n}\binom{n}{k}B_{k}(w_{1}y_{1})B_{n-k}(y_{2})w_{1}^{n-k}\\
&=\sum_{k=0}^{n}\binom{n}{k}B_{k}(y_{1})B_{n-k}(w_{1}y_{2})w_{1}^{k}\\
\end{split}
\end{equation*}
\begin{equation*}
=\sum_{k+\ell+m=n}\binom{n}{k,\ell,m}B_{k}(y_{1})B_{\ell}(y_{2})S_{m}(w_{1}-1)w_{1}^{k+\ell-1}.
\end{equation*}

{\THM {Theorem {5}}}{\it Let $w_{1}, w_{2},w_{3} $ be any positive integers.
Then the following expression  is invariant under any permutation of
$w_{1}, w_{2},w_{3}$, so that it gives us six symmetries.
\rm}

\begin{equation*}
\tag{48}
\begin{split}
& w_{1}^{n-1} \sum_{k=0}^{n}\binom{n}{k}B_{k}(w_{3}y_{1})\sum_{i=0}^{w_{1}-1}B_{n-k}(w_{2}y_{2}+\frac{w_{2}}{w_{1}}i)
w_{3}^{n-k}w_{2}^{k}\\
&=w_{1}^{n-1} \sum_{k=0}^{n}\binom{n}{k}B_{k}(w_{2}y_{1})\sum_{i=0}^{w_{1}-1}B_{n-k}(w_{3}y_{2}+\frac{w_{3}}{w_{1}}i)
w_{2}^{n-k}w_{3}^{k} \\
&=w_{2}^{n-1} \sum_{k=0}^{n}\binom{n}{k}B_{k}(w_{3}y_{1})\sum_{i=0}^{w_{2}-1}B_{n-k}(w_{1}y_{2}+\frac{w_{1}}{w_{2}}i)
w_{3}^{n-k}w_{1}^{k}\\
&=w_{2}^{n-1} \sum_{k=0}^{n}\binom{n}{k}B_{k}(w_{1}y_{1})\sum_{i=0}^{w_{2}-1}B_{n-k}(w_{3}y_{2}+\frac{w_{3}}{w_{2}}i)
w_{1}^{n-k}w_{3}^{k} \\
&=w_{3}^{n-1} \sum_{k=0}^{n}\binom{n}{k}B_{k}(w_{2}y_{1})\sum_{i=0}^{w_{3}-1}B_{n-k}(w_{1}y_{2}+\frac{w_{1}}{w_{3}}i)
w_{2}^{n-k}w_{1}^{k}\\
\end{split}
\end{equation*}
\begin{equation*}
=w_{3}^{n-1} \sum_{k=0}^{n}\binom{n}{k}B_{k}(w_{1}y_{1})\sum_{i=0}^{w_{3}-1}B_{n-k}(w_{2}y_{2}+\frac{w_{2}}{w_{3}}i)
w_{1}^{n-k}w_{2}^{k}.
\end{equation*}

Letting  $w_{3}=1$ in (48), we obtain alternative expressions for the identities in (46).

{\COR {Corollary {6}}} {\it Let $w_{1}, w_{2}$ be any positive integers. \rm}
\begin{equation*}
\tag{49}
\begin{split}
&\sum_{k=0}^{n}\binom{n}{k}B_{k}(w_{1}y_{1})B_{n-k}(w_{2}y_{2})w_{1}^{n-k}w_{2}^{k}\\
&=\sum_{k=0}^{n}\binom{n}{k}B_{k}(w_{2}y_{1})B_{n-k}(w_{1}y_{2})w_{2}^{n-k}w_{1}^{k}\\
&=w_{1}^{n-1}\sum_{k=0}^{n}\binom{n}{k}B_{k}(y_{1})\sum_{i=0}^{w_{1}-1}B_{n-k}(w_{2}y_{2}+\frac{w_{2}}{w_{1}}i)w_{2}^{k}\\
\end{split}
\end{equation*}

\begin{equation*}
\begin{split}
&=w_{1}^{n-1}\sum_{k=0}^{n}\binom{n}{k}B_{k}(w_{2}y_{1})\sum_{i=0}^{w_{1}-1}B_{n-k}(y_{2}+\frac{i}{w_{1}})w_{2}^{n-k}\\
&=w_{2}^{n-1}\sum_{k=0}^{n}\binom{n}{k}B_{k}(y_{1})\sum_{i=0}^{w_{2}-1}B_{n-k}(w_{1}y_{2}+\frac{w_{1}}{w_{2}}i)w_{1}^{k}\\
&=w_{2}^{n-1}\sum_{k=0}^{n}\binom{n}{k}B_{k}(w_{1}y_{1})\sum_{i=0}^{w_{2}-1}B_{n-k}(y_{2}+\frac{i}{w_{2}})w_{1}^{n-k}.\\
\end{split}
\end{equation*}

Putting further $w_{2}=1$ in (49), we have the alternative expressions for the identities for (47).

{\COR {Corollary {7}}} {\it Let $w_{1}$ be any positive integer. \rm}

\begin{equation*}
\begin{split}
&\sum_{k=0}^{n}\binom{n}{k}B_{k}(y_{1})B_{n-k}(w_{1}y_{2})w_{1}^{k}\\
&=\sum_{k=0}^{n}\binom{n}{k}B_{k}(y_{2})B_{n-k}(w_{1}y_{1})w_{1}^{k}\\
&=w_{1}^{n-1}\sum_{k=0}^{n}\binom{n}{k}B_{k}(y_{1})\sum_{i=0}^{w_{1}-1}B_{n-k}(y_{2}+\frac{i}{w_{1}}).\\
\end{split}
\end{equation*}

{\THM {Theorem {8}}}{\it Let $w_{1}, w_{2},w_{3} $ be any positive integers.
Then we have the following three symmetries in $w_{1}, w_{2},w_{3}$:
\rm}
\begin{equation*}
\tag{50}
\sum_{k+\ell+m=n}\binom{n}{k,\ell,m}B_{k}(w_{1}y_{1})S_{\ell}(w_{2}-1)S_{m}(w_{3}-1)
w_{1}^{\ell+m}w_{2}^{k+m-1}w_{3}^{k+\ell -1}
\end{equation*}
\begin{equation*}
\tag{51}
=\sum_{k+\ell+m=n}\binom{n}{k,\ell,m}B_{k}(w_{2}y_{1})S_{\ell}(w_{3}-1)S_{m}(w_{1}-1)
w_{2}^{\ell+m}w_{3}^{k+m-1}w_{1}^{k+\ell -1}
\end{equation*}
\begin{equation*}
\tag{52}
=\sum_{k+\ell+m=n}\binom{n}{k,\ell,m}B_{k}(w_{3}y_{1})S_{\ell}(w_{1}-1)S_{m}(w_{2}-1)
w_{3}^{\ell+m}w_{1}^{k+m-1}w_{2}^{k+\ell -1}.
\end{equation*}

Putting $w_{3}=1$ in (50)-(52), we get the following corollary.

{\COR {Corollary {9}}} {\it Let $w_{1}, w_{2}$ be any positive integers. \rm}
\begin{equation*}
\tag{53}
\begin{split}
&\sum_{k=0}^{n}\binom{n}{k}B_{k}(w_{1}y_{1})S_{n-k}(w_{2}-1)w_{1}^{n-k}w_{2}^{k-1}\\
&=\sum_{k=0}^{n}\binom{n}{k}B_{k}(w_{2}y_{1})S_{n-k}(w_{1}-1)w_{2}^{n-k}w_{1}^{k-1}\\
&=\sum_{k+\ell+m=n}\binom{n}{k,\ell,m}B_{k}(y_{1})S_{\ell}(w_{1}-1)S_{m}(w_{2}-1)w_{1}^{k+m-1}w_{2}^{k+\ell-1}.\\
\end{split}
\end{equation*}

Letting further $w_{2}=1$ in (53), we get the following corollary.
This is also obtained in [6, Cor.2] and mentioned in [3].

{\COR {Corollary {10}}} {\it Let $w_{1}$ be any positive integer. \rm}

\begin{equation*}
\tag{54}
B_{n}(w _{1}y_{1})=\sum_{k=0}^{n}\binom{n}{k}B_{k}(y_{1})S_{n-k}(w_{1}-1)w_{1}^{k-1}.\\
\end{equation*}

{\THM {Theorem {11}}}{\it Let $w_{1}, w_{2},w_{3} $ be any positive integers.
Then the following expression  is invariant under any permutation of
$w_{1}, w_{2},w_{3}$, so that it gives us six symmetries.
\rm}

\begin{equation*}
\tag{55}
\begin{split}
& w_{1}^{n-1}\sum_{k=0}^{n}\binom{n}{k}\sum_{i=0}^{w_{1}-1}B_{k}(w_{2}y_{1}+\frac{w_{2}}
{w_{1}}i)S_{n-k}(w_{3}-1)w_{2}^{n-k}w_{3}^{k-1}\\
&=w_{1}^{n-1}\sum_{k=0}^{n}\binom{n}{k}\sum_{i=0}^{w_{1}-1}B_{k}(w_{3}y_{1}+\frac{w_{3}}
{w_{1}}i)S_{n-k}(w_{2}-1)w_{3}^{n-k}w_{2}^{k-1}\\
&=w_{2}^{n-1}\sum_{k=0}^{n}\binom{n}{k}\sum_{i=0}^{w_{2}-1}B_{k}(w_{1}y_{1}+\frac{w_{1}}{w_{2}}i)S_{n-k}(w_{3}-1)
w_{1}^{n-k}w_{3}^{k-1}\\
&=w_{2}^{n-1}\sum_{k=0}^{n}\binom{n}{k}\sum_{i=0}^{w_{2}-1}B_{k}(w_{3}y_{1}+\frac{w_{3}}{w_{2}}i)S_{n-k}(w_{1}-1)
w_{3}^{n-k}w_{1}^{k-1}\\
&=w_{3}^{n-1}\sum_{k=0}^{n}\binom{n}{k}\sum_{i=0}^{w_{3}-1}B_{k}(w_{1}y_{1}+\frac{w_{1}}{w_{3}}i)S_{n-k}(w_{2}-1)
w_{1}^{n-k}w_{2}^{k-1}\\
\end{split}
\end{equation*}

\begin{equation*}
\begin{split}
&=w_{3}^{n-1}\sum_{k=0}^{n}\binom{n}{k}\sum_{i=0}^{w_{3}-1}B_{k}(w_{2}y_{1}+\frac{w_{2}}{w_{3}}i)S_{n-k}(w_{1}-1)
w_{2}^{n-k}w_{1}^{k-1}.\\
\end{split}
\end{equation*}

Putting $w_{3}=1$ in (55), we obtain the following corollary. In Section I, the identities in (53), (56), and (58) are combined to give those in (8)-(15).

{\COR {Corollary {12}}} {\it Let $w_{1}, w_{2}$ be any positive integers. \rm}
\begin{equation*}
\tag{56}
\begin{split}
&w_{1}^{n-1}\sum_{i=0}^{w_{1}-1}B_{n}(w_{2}y_{1}+\frac{w_{2}}{w_{1}}i)\\
&=w_{2}^{n-1}\sum_{i=0}^{w_{2}-1}B_{n}(w_{1}y_{1}+\frac{w_{1}}{w_{2}}i)\\
&=\sum_{k=0}^{n}\binom{n}{k}B_{k}(w_{2}y_{1})S_{n-k}(w_{1}-1)w_{2}^{n-k}w_{1}^{k-1}\\
&=\sum_{k=0}^{n}\binom{n}{k}B_{k}(w_{1}y_{1})S_{n-k}(w_{2}-1)w_{1}^{n-k}w_{2}^{k-1}\\
&=w_{1}^{n-1}\sum_{k=0}^{n}\binom{n}{k}\sum_{i=0}^{w_{1}-1}B_{k}(y_{1}+\frac{i}{w_{1}})S_{n-k}(w_{2}-1)w_{2}^{k-1}\\
\end{split}
\end{equation*}
\begin{equation*}
=w_{2}^{n-1}\sum_{k=0}^{n}\binom{n}{k}\sum_{i=0}^{w_{2}-1}B_{k}(y_{1}+\frac{i}{w_{2}})S_{n-k}(w_{1}-1)w_{1}^{k-1}.
\end{equation*}

Letting further $w_{2}=1$ in (56), we get the following corollary. This is the well-known multiplication formula for Bernoulli polynomials together with the relatively new identity mentioned in (54).

{\COR {Corollary {13}}} {\it Let $w_{1}$ be any positive integer. \rm}
\begin{equation*}
\begin{split}
B_{n}(w _{1}y_{1})
&= w_{1}^{n-1}\sum_{i=0}^{w_{1}-1}B_{n}(y_{1}+\frac{i}{w_{1}})\\
&= \sum_{k=0}^{n}\binom{n}{k}B_{k}(y_{1})S_{n-k}(w_{1}-1)w_{1}^{k-1}.\\
\end{split}
\end{equation*}

{\THM {Theorem {14}}}{\it Let $w_{1}, w_{2},w_{3} $ be any positive integers.
Then we have the following three symmetries in $w_{1}, w_{2},w_{3}$ :
\rm}
\begin{equation*}
\tag{57}
\begin{split}
& (w_{1}w_{2})^{n-1}\sum_{i=0}^{w_{1}-1}\sum_{j=0}^{w_{2}-1}B_{n}(w_{3}y_{1}+\frac{w_{3}}
{w_{1}}i+\frac{w_{3}}{w_{2}}j)\\
&=(w_{2}w_{3})^{n-1}\sum_{i=0}^{w_{2}-1}\sum_{j=0}^{w_{3}-1}B_{n}(w_{1}y_{1}+\frac{w_{1}}
{w_{2}}i+\frac{w_{1}}{w_{3}}j)\\
&=(w_{3}w_{1})^{n-1}\sum_{i=0}^{w_{3}-1}\sum_{j=0}^{w_{1}-1}B_{n}(w_{2}y_{1}+\frac{w_{2}}
{w_{3}}i+\frac{w_{2}}{w_{3}}j).\\
\end{split}
\end{equation*}

Letting $w_{3}=1$ in (57), we have the following corollary.

{\COR {Corollary {15}}} {\it Let $w_{1}, w_{2}$ be any positive integers. \rm}
\begin{equation*}
\tag{58}
\begin{split}
&w_{1}^{n-1}\sum_{j=0}^{w_{1}-1}B_{n}(w _{2}y_{1}+\frac{w_{2}}{w_{1}}j)\\
&= w_{2}^{n-1}\sum_{i=0}^{w_{2}-1}B_{n}(w_{1}y_{1}+\frac{w_{1}}{w_{2}}i)\\
&=(w_{1}w_{2})^{n-1}\sum_{i=0}^{w_{1}-1}\sum_{j=0}^{w_{2}-1}B_{n}(y_{1}+\frac{i}{w_{1}}+\frac{j}{w_{2}}).\\
\end{split}
\end{equation*}

{\THM {Theorem {16}}}{\it Let $w_{1}, w_{2},w_{3} $ be any positive integers.
Then we have the following two symmetries in $w_{1}, w_{2},w_{3}$:
\rm}
\begin{equation*}
\tag{59}
\begin{split}
&\sum_{k+\ell+m=n}\binom{n}{k,\ell,m}B_{k}(w _{1}y)B_{\ell}(w _{2}y)B_{m}(w _{3}y)w_{3}^{k}w_{1}^{\ell}w_{2}^{m}\\
&= \sum_{k+\ell+m=n}\binom{n}{k,\ell,m}B_{k}(w _{1}y)B_{\ell}(w _{3}y)B_{m}(w _{2}y)w_{2}^{k}w_{1}^{\ell}w_{3}^{m}.\\
\end{split}
\end{equation*}

{\THM {Theorem {17}}}{\it Let $w_{1}, w_{2},w_{3} $ be any positive integers.
Then we have the following two symmetries in $w_{1}, w_{2},w_{3}$:
\begin{equation*}
\tag{60}
\sum_{k+\ell+m=n}\binom{n}{k,\ell,m}S_{k}(w_{1}-1)S_{\ell}(w_{2}-1)S_{m}(w_{3}-1)
w_{3}^{k-1}w_{1}^{\ell-1}w_{2}^{m-1}
\end{equation*}
\begin{equation*}
\tag{61}
= \sum_{k+\ell+m=n}\binom{n}{k,\ell,m}S_{k}(w _{1}-1)S_{\ell}(w _{3}-1)S_{m}(w _{2}-1)w_{2}^{k-1}w_{1}^{\ell-1}w_{3}^{m-1}.
\end{equation*}
\rm}

Putting $w_{3}=1$ in (60) and (61) and multiplying the resulting identity by $w_{1} w_{2}$, we get the following corollary.

{\COR {Corollary {18}}} {\it Let $w_{1}, w_{2}$ be any positive integers. \rm}
\begin{equation*}
\begin{split}
&\sum_{k=0}^{n}\binom{n}{k}S_{k}(w _{2}-1)S_{n-k}(w_{1}-1)w_{1}^{k}\\
&=\sum_{k=0}^{n}\binom{n}{k}S_{k}(w _{1}-1)S_{n-k}(w_{2}-1)w_{2}^{k}. \\
\end{split}
\end{equation*}

\bigskip\bigskip
\begin{large}
{\sc References}
\end{large}
\par
\begin{enumerate}
\item[{[1]}] E. Deeba and D. Rodriguez, {\it Stirling's and Bernoulli numbers},
Amer. Math. Monthly 98(1991), 423-426.

\item[{[2]}] F. T. Howard, {\it Applications of a recurrence for the Bernoulli numbers},
J. Number Theory 52(1995), 157-172.

\item[{[3]}] T. Kim, {\it Symmetry $p$-adic invariant integral on $\mathbb{Z}_{p}$ for Bernoulli and Euler polynomials}, J. Adv. Difference Equ. Appl. 14(2008), 1267-1277.

\item[{[4]}] W. H. Schikhof, {\it "Ultrametric calculus: An introduction to $p$-adic analysis,"} Cambridge University Press, 2006.

\item[{[5]}] H. Tuenter, {\it A symmetry of power sum polynomials and Bernoulli numbers},
Amer. Math. Monthly 108(2001), 258-261.

\item[{[6]}] S. Yang, {\it An identity of symmetry for the Bernoulli polynomials},
Discrete Math. 308(2008), 550-554.

\end{enumerate}

\end{document}